\documentclass[12pt,draftcls,onecolumn]{IEEEtran}

\ifCLASSOPTIONcompsoc
 
  \usepackage[nocompress]{cite}
\else
  \usepackage{cite}
\fi

\ifCLASSINFOpdf
  
\else

\fi

\usepackage{algorithm}
 \usepackage{algpseudocode}
 \usepackage{amsmath}
 \usepackage{graphics}
 \usepackage{epsfig}
 \usepackage[misc]{ifsym}
 \usepackage{bbding}

\usepackage{graphicx}
\usepackage{epsf}
\usepackage{amssymb}

\newtheorem{thm}{Theorem}[section]

\newtheorem{cor}[thm]{Corollary}
\newtheorem{lem}[thm]{Lemma}
\newtheorem{prop}[thm]{Proposition}

\newtheorem{HDRthm}[thm]{High Dimension Roll Theorem}
\newtheorem{Ethm}[thm]{Essen Theorem}

\newtheorem{PAthm}[thm]{Polynomial Automorphism Theorem}

\newtheorem{HTthm}[thm]{Hadamard Theorem}
\newtheorem{Pthm}[thm]{Proper Map Theorem}

\begin{document}

\title{An Optimization Approach to Jacobian Conjecture}

\author{\IEEEauthorblockN{Jiang Liu}\\
\IEEEauthorblockA{Chongqing Institute of Green and Intelligent Technology\\
Chinese Academy of Sciences\\
\Letter: liujiang@cigit.ac.cn}
}

\maketitle

\begin{abstract}
Let $n\geq 2$ and $\mathbb K $ be a number field of characteristic $0$. Jacobian Conjecture asserts for a polynomial map $\mathcal P$ from $\mathbb K ^n$ to itself, if the determinant of its Jacobian matrix is a nonzero constant in $\mathbb K $ then the inverse $\mathcal P^{-1}$ exists and is also a polynomial map. This conjecture was firstly proposed by Keller in 1939 for $\mathbb K ^n=\mathbb C^2$ and put in Smale's 1998 list of Mathematical Problems for the Next Century. This study is going to present a proof for  the conjecture. Our proof is based on Dru{\.{z}}kowski Map and Hadamard's Diffeomorphism Theorem, and additionally uses some optimization idea.
\end{abstract}


\IEEEpeerreviewmaketitle

\section{Introduction}\label{sec:intro}
Let  $\mathbb K $ denote a number field, on which  $\mathbb K [X] =\mathbb K [x_1,\dots, x_n]$ be the polynomial ring in $n$ variables $X=(x_1,\dots, x_n)$.  Each polynomial vector $\mathcal P(X):=(\mathcal P_1(X),\dots, \mathcal P_n(X))\in  \mathbb K [X]^n$ defines a map from $\mathbb K ^n$ to  $\mathbb K ^n$. The $n\times n$ Jacobian matrix $\mathbf J_{\mathcal P}(X)$ ($\mathbf{J}_{\mathcal P}$ for short) of $\mathcal P(X)$ consists of the partial derivatives of $\mathcal P_j$ with respect to $x_k$, that is, $\mathbf{J}_{\mathbf P}=[\frac{\partial  P_j}{\partial  x_k},j,k=1,\dots,n]$.
Then the determinant $det(\mathbf{J}_{\mathcal P})$ of $\mathbf{J}_{\mathcal P}$ is a polynomial function of $X$. \textit{Jacobian Condition} is that $det(\mathbf{J}_{\mathcal P})$ is a nonzero constant in $\mathbb K $. A $\mathcal P\in \mathbb K [X]^n$ is called \textit{Keller map} \cite{Essenbook2000} if it satisfies Jacobian Condition. A polynomial map $\mathcal P\in \mathbb K [X]^n$  is an \textit{automorphism} of $\mathbb K ^n$ if the inverse $\mathcal P^{-1}$  exists and is also a polynomial map. Due to Osgood's Theorem, the Jacobian Condition is necessary for  $\mathcal P$ being an automorphism.  Keller  \cite{Keller1939}  proposed the Jacobian Conjecture ($\mathcal{JC}$) for $\mathbb K^n =\mathbb C^2$  in 1939,  where $\mathbb C$ is the complexity field.  We refer the reader to \cite{Bass1982formalinverse,Druzkowski1997ReviewJC,Essenbook2000,denessen2001review,Wright2005Review} for a nice survey paper containing some history and updated progresses  of the Jacobian Conjecture. Abhyankar \cite{Abhyankar1977} gave its modern style as follows.

{\textit{Jacobian Conjecture} --- $\mathcal{JC}(\mathbb K , n)$:} If $\mathcal P\in \mathbb K [X]^n$ is a Keller map where $\mathbb K $ is an arbitrary field of characteristic $0$ and the integer $n\geq 2$, then  $\mathcal P$ is an automorphism.

$\mathcal{JC}(\mathbb K, 1)$ is trivially true.
Besides, there are simple counterexamples  \cite{Abhyankar1974Inversion,Bass1982formalinverse} for $\mathcal{JC}(\mathbb K , n) $  when the number field  $\mathbb K $ has characteristic $>0$. So the characteristic $0$ condition is necessary. These two assumptions are always supposed to be true in the rest of this article unless a related statement is specified particularly. Let $\mathcal{JC}(\mathbb K)$ stand for $\mathcal{JC}(\mathbb K , n)$ being true for all positive integers $n$. Due to Lefschetz Principle, it suffices to deal only with $\mathcal{JC}(\mathbb C)$. For polynomial maps in $\mathbb C [X]^n$, it is well known that
\begin{thm}[\cite{Winiarski1979Inverse,Yagzhev1981,Bass1982formalinverse,RudinInjective1995}]
	If a polynomial map $\mathbf P\in \mathbb C[X]^n$ is injective then it must be an automorphism.	
\end{thm}

Thus, it suffices to show that each Keller map in $\mathbb C[X]^n$ is injective. There is a remarkable result \cite{Wang1980degree2} for polynomial maps of degree $2$.
\begin{thm}\label{thm:2deg}
	If $\mathcal P(X)\in \mathbb K[X]^n$ is a quadratic polynomial, then $\mathcal P(X)$ is injective iff for each $X\in \mathbb K$ it alway has $det(\mathbf J_{\mathcal P}(X))\neq 0$.
\end{thm}
\begin{IEEEproof}
For the result, we present the wonderful proof gave in \cite{Ludwik2005}. The assertion is true by
\begin{equation}
\mathcal P(X)-\mathcal P(Y)=\mathbf J_{\mathcal P}\left(\frac{X+Y}{2}\right)(X-Y)
\end{equation}
		
\end{IEEEproof}	

For high degree polynomial maps, Yag{\v{z}}ev \cite{Yagzhev1981} and independently Bass et al. \cite{Bass1982formalinverse} showed that it is sufficient to consider Keller maps of form $\mathcal P(X):=X+\mathcal H(X)\in \mathbb C[X]^n$ which is called \textit{Yag{\v{z}}ev map}, where $\mathcal H(X)$ is a homogeneous polynomial vectors of degree $3$.  Based on this result,  Dru{\.{z}}kowski \cite{druzkowski1983} further  showed  that it suffices to consider  Keller maps of so-called cubic linear form $\mathcal P(X):=X+(\mathcal AX)^{*3}\in \mathbb C[X]^n$, nowadays, called {\textit{Dru{\.{z}}kowski map}} (D-map for short),  wherein $(\mathcal AX)^{*3}$ be the vector whose $k$-th element is $(A_k^TX)^3$ and $A^T_k$ is the $k$-th row vector of $\mathcal A$. That is, if each D-map $\mathcal P(X)\in \mathbb C[X]^n$ is injective then Jacobian Conjecture is true. It is natural to ask whether all D-map $\mathcal P(X)\in \mathbb R[X]^n$ being injective also implies Jacobian Conjecture. The answer is yes in what follows.

\begin{thm}\label{thm:inj}
	If each  D-map in $\mathbb R[X]^n$ is injective for all $n$, then each Keller map $\mathbf P\in \mathbb C[X]^n$ is injective and so an automorphism. 
\end{thm}
\begin{IEEEproof}
Assume variable $X:=Y+\mathbf i Z$ with $\mathbf i=\sqrt{-1}$, that is, $x_k=y_k+ \mathbf i z_k$ for $1\leq k\leq n$. Like the argument on the first page in \cite{druzkowski1995SurveyCounterexamples} or in the proof of Lemma 6.3.13 on page 130 in \cite{Essenbook2000}, a D-map  $\mathbf P(X):= \mathbf P(Y+ \mathbf i Z)\in\mathbb C[X]^n$ can be treated as a real map $\mathbf P^\star(Y,Z):=(Re P_1(Y+ \mathbf iZ), Im  P_1(Y+\mathbf iZ), \cdots, Re P_n(Y+\mathbf iZ), Im  P_n(Y+ \mathbf iZ))\in \mathbb R[Y,Z]^{2n}$, where $Re P_k(Y+\mathbf iZ)$ and $ Im  P_k(Y+\mathbf iZ)$ are the real and imaginary  parts of $P_k(Y+\mathbf iZ)$, respectively. It is evident that $\mathbf P(Y+\mathbf iZ)\in\mathbb C[X]^n$ and $\mathbf P^\star(Y,Z)\in \mathbb R[Y,Z]^{2n}$ have the same injectivity. In fact, each $Re P_k(Y+\mathbf iZ)$ and $Im  P_k(Y+\mathbf iZ)$ for $1\leq k\leq n$ are of forms $Y_k+H_k (Y,Z)$ and $Z_k+H_k^*(Y,Z)$, respectively. Furthermore, both $H_k(Y,Z)$ and $H_k^*(Y,Z)$ are homogeneous polynomial about variables $(Y,Z)$ of degree $3$, which can be seen by expanding the $k$-th component of  D-map $P_k(Y+\mathbf iZ)$. It is well known $det(\mathbf J_{\mathbf P^\star}(Y,Z))=|det(\mathbf J_{\mathbf P}(Y+\mathbf iZ))|^2$. As $\mathbf P(Y+\mathbf iZ)$ is a D-map, $det(\mathbf J_{\mathbf P^\star}(Y,Z))=1$ and thus $\mathbf P^\star(Y,Z)$ is actually a Yag{\v{z}}ev map. Now, we apply Dru{\.{z}}kowski's reduction in \cite{druzkowski1983} to get that there is some D-map $\mathbf P^{@}(Y,Z,\tilde{Z})\in \mathbb R[Y,Z,\tilde{Z}]^{m}$ with $m\geq n$ such that  $\mathbf P^\star(Y,Z)$ is injective iff  $\mathbf P^{@}(Y,Z,\tilde{Z})$ is injective,  by Theorem 3 and Remark 4 in \cite{druzkowski1983}. In consequence, if each  D-map in $\mathbb R[X]^n$ is injective for all $n$, then each Keller map $\mathbf P\in \mathbb C[X]^n$ is injective and so an automorphism.
\end{IEEEproof}

As a result of Theorem \ref{thm:inj}, Jacobian Conjecture is true if each  D-map in $\mathbb R[X]^n$ is injective for all $n$. In the study, we will use Hadamard's Diffeomorphism Theorem to show the injectivity of D-maps in $\mathbb R[X]^n$. By Jacobian Condition, each Keller map is a local homeomorphism. Even for analytic maps,
Hadamard's Diffeomorphism Theorem presents a necessary and sufficient condition for local homeomorphism being diffeomorphisme in terms of the proper map. A differentiable map $\mathcal F: \mathbb R^n\mapsto \mathbb R^n$ is an \textit{diffeomorphism} if it is bijective and its inverse is also differentiable. A continuous map $\mathcal H: \mathbb R^n\mapsto \mathbb R^n$ is \textit{proper} if $\mathcal H^{-1}(\mathfrak C)$ is compact for each compact set $\mathfrak C\subset \mathbb R^n$. Equivalently, a map $\mathcal H$ is proper iff it maps each unbounded set $\mathfrak B$ into an unbounded set $\mathcal H(\mathfrak B)$. Based on these notions, in (\cite{Hadamard1906,Gordon1972On}) it was showed that
\begin{HTthm}\label{thm:HTthm}
For an analytic map $\mathcal H: \mathbb R^n\mapsto \mathbb R^n$, $\mathcal H$ is a diffeomorphism iff $\det\mathbf J_{\mathcal H}(X)\neq 0$ for all  $X\in \mathbb R^n$ and  $\mathcal H$ is proper.
\end{HTthm}

Therefore, it suffices to show the properness of D-maps in $\mathbb R[X]^n$ for their injectivity. To this end, we are going to show the following crucial result in the next section.

\begin{Pthm}\label{thm:D-proper}
Each D-map $\mathcal D(X):=X+(\mathcal A X)^{*3}$ in $\mathbb R[X]^n$ is proper and so injective.
\end{Pthm}

The basic idea of its proof is in what follows. Firstly, we show $\mathcal I+ \lambda\mathcal A(\mathcal AX)^{\Delta2}$ is invertible for each $\lambda\in \mathbb R$, where $\mathcal I$ is the identity matrix. Then we show the result by contradiction. Roughly, if  $\mathcal D(X)$ is not proper, then we will obtain a point $Z=\mathcal A Y$ such that $Z+\ell \mathcal A Z^{*3}=\mathbf 0$ for some $\ell\in\mathbb R$ and $Y\in\mathbb R^n$. The detail will be depicted in the next section. By Theorem \ref{thm:inj} we infer

\begin{PAthm}\label{thm:Conjecture}
For every $n\geq 2$ and any field $ \mathbb K $ of characteristic $0$, each Keller map $\mathcal P\in \mathbb K [X]^n$ is an automorphism.	
\end{PAthm}

In the rest of this article, we assume $n\geq 2$ for the dimension of $\mathbb R^n$ unless it is specified in particular. For clarity, let capital letters denote vectors, blackboard bold letters denote sets of vectors (points), fraktur letters denote vector functions, little letters denote integers, Greek letters denote real numbers and real functions.

\section{Proof of Proper Map Theorem}

Let's start with some basic notions and notations.  A set $\mathfrak B\subseteq \mathbb R^n$ is called {\textit{unbunded}} if for any $\gamma\in\mathbb R$ there is an element $X\in \mathfrak B$ such that $||X||>\gamma$, where $||\cdot||$ is the Euclidean norm, that is, $||X||:=\sqrt{\sum\limits_{1\le k\leq n}X_k^2}$. Let  $\mathfrak U:=\{X\in\mathbb R^n\mid ||X||=1\}$, which is a compact set in $\mathbb R^n$. Given a matrix $\mathcal A$, let $\mathcal A^T$ denote its transpose, define $Im (\mathcal A):=\{\mathcal AX \mid X\in\mathbb R^n\}$. Given a vector $V=(v_1,\cdots, v_n)^T$, let $V^{\Delta2}$ denote the diagonal matrix whose diagonals are $v_1^2, \dots, v_n^2$. For a D-map $\mathcal D(X):=X+(\mathcal AX)^{*3}$, its Jacobian matrix is $\mathbf J_{\mathcal P}=\mathcal I+3(\mathcal AX)^{\Delta2}\mathcal A$. From \cite{druzkowski1983},

\begin{thm}\label{thm:Nilpthm}
	$\mathcal D(X):=X+(\mathcal AX)^{*3}$ is a D-map iff $(\mathcal AX)^{\Delta2}\mathcal A$ is  nilpotent matrix for any $X$. 
\end{thm}

From this theorem, the following proposition is evident.
\begin{prop}\label{prop:d-mapclas}
	If  $\mathcal D(X):=X+(\mathcal AX)^{*3}$ is a D-map, then for all $\lambda\in \mathbb R$, the map $\mathcal  D_\lambda(X):=X+\lambda(\mathcal AX)^{*3}$ is also a D-map. Thus, $\lambda(\mathcal AX)^{\Delta2}\mathcal A$ and so $\lambda\mathcal A(\mathcal AX)^{\Delta2}$ are nilpotent matrices. Therefore, $\mathcal I +\lambda\mathcal A(\mathcal AX)^{\Delta2}$ is invertible. 
\end{prop}
\begin{IEEEproof}
	By Theorem \ref{thm:Nilpthm}, if  $\mathcal D(X):=X+(\mathcal AX)^{*3}$ is a D-map, then $(\mathcal AX)^{\Delta2}\mathcal A$ is nilpotent matrix,  and so is $\lambda(\mathcal AX)^{\Delta2}\mathcal A$ for any  $\lambda\in \mathbb R$. It is evident that $(\mathcal A\mathcal B)^k=0_{n\times n}$ implies $(\mathcal B\mathcal A)^{k+1}=\mathcal B (\mathcal A\mathcal B)^k\mathcal A=0_{n\times n}$ for any integer $k$ and arbitrary matrices $\mathcal A,\mathcal B\in\mathbb R^{n\times n}$, where $0_{n\times n}$ is the zero matrix of size $n\times n$. Hence, all $\lambda\mathcal A(\mathcal AX)^{\Delta2}$ are nilpotent matrices. In consequence, each $\mathcal I +\lambda\mathcal A(\mathcal AX)^{\Delta2}$ is invertible.
	
	\end{IEEEproof}

This result accomplished the first part of the proof for Theorem \ref{thm:D-proper}. To finish the proof,  we are going to show

\begin{lem}\label{lem:nonprop}
If  a linear cubic form  $\hat{\mathcal C}(X):=X+(\mathcal AX)^{*3}$ is not proper, then (I) there is an unbounded $\{U_i\in Im(\mathcal A\mathcal A^T)\}_i$ such that $\{U_i+\mathcal A U_i^{*3}\}_i$ is bounded, that is, $\hat{\mathcal C}(X)$ is not proper on $Im(\mathcal A\mathcal A^T)$; and (II) there is some $W\in Im(\mathcal A^T)$ such that $\mathcal A(\mathcal AW)^{*3}=0_{n\times n}$ but $\mathcal AW\neq 0_{n\times n}$.
\end{lem}

To show this result, we need following preliminary materials. Given a matrix $\mathcal A\in \mathbb R^{n\times n}$, let $\mathfrak  A^\bot:=\{X\in \mathbb R^n\mid \mathcal A X=\mathbf 0\}$ denote the linear space in which each element is a solution of $\mathcal A X=\mathbf 0$. For any $ X\in  \mathfrak A^\bot$ and $Y=\mathcal A^TZ\in Im(\mathcal A^T)$, we have $Y^TX=Z^T\mathcal A X=0$, that is, $Y\bot X$. In consequence,
\begin{prop}\label{prop:property} The spaces $ \mathfrak A^\bot$ and $Im(\mathcal A^T)$ have the following properties:
	\begin{enumerate}
		\item $\mathfrak A^\bot$ and $Im(\mathcal A^T)$ are orthogonal to each other.
		\item $\mathfrak A^\bot \cap Im(\mathcal A^T) =\mathbf 0$, and so $|| \mathcal AY||\neq 0$ for any nonzero $Y\in Im(\mathcal A^T) $.
		\item Both $\mathfrak A^\bot$ and $Im(\mathcal A^T)$ are closed sets.
		\item For any nonzero $Z\in \mathbb R^n$, it can be uniquely decomposed as $Z=X+Y$ such that $X\bot Y$, $X\in \mathfrak A^\bot$ and $Y\in Im(\mathcal A^T)$.
	\end{enumerate}
\end{prop}

Let's set  $\mathfrak U(\mathfrak A^\bot):=\mathfrak U \cap \mathfrak A^\bot$ and $\mathfrak U(Im(\mathcal A^T)):=\mathfrak U \cap Im(\mathcal A^T)$, then
\begin{prop}\label{prop:unitclosed}
	Both $\mathfrak U(\mathfrak A^\bot)$ and $\mathfrak U(Im(\mathcal A^T))$ are closed  and so compact. Moreover, $||\mathcal AY||\neq 0$ for any $Y\in \mathfrak  U(Im(\mathcal A^T))$.
\end{prop}	

Now, we present the proof of Lemma \ref{lem:nonprop} in what follows.

\begin{IEEEproof}
	If $\hat{\mathcal C}(X)$ is not proper, then there  must  be an unbounded sequence $\{Y_i\}_i$  such that $\{\mathcal C(Y_i)\}_i$ is a bounded sequence.  Without loss of generality, we suppose that $||Y_{i+1}||>||Y_i||$ for all $i$. By the boundness of $\{\mathcal C(Y_i)\}_i$, there is some real number $\sigma<\infty$ such that
	\begin{equation}\label{ineq:bund}
	Y_i^TY_i+2Y_i^T(\mathcal AY_i)^{*3}+||(\mathcal AY_i)^{*3}||^2\leq \sigma^2
	\end{equation}
	for all $i$. Because $\{Y_i\}_i$ is unbounded and and $\{\mathcal C(Y_i)\}_i$ is  bounded, it must be $Y_i^T(\mathcal AY_i)^{*3}<0$ for almost all $i$.  Without loss of generality, we suppose  $Y_i^T(\mathcal AY_i)^{*3}<0$   for all $i$. As a result, $(\mathcal  AY_i)^{*3}\neq \mathbf 0$ and so $||(\mathcal  AY_i)^{*3}||>0$ for all $i$. Then by optimization theory on quadratic polynomials \cite{Boyd2004Optimization},  we have
	\begin{equation}
	||(\mathcal  AY_i)^{*3}||^2\gamma^2+2Y_i^T(\mathcal AY_i)^{*3}\gamma+Y_i^TY_i \geq  Y_i^TY_i-\frac{(Y_i^T(\mathcal AY_i)^{*3})^2}{||(\mathcal AY_i)^{*3}||^2}
	\end{equation}
	for any $\gamma\in \mathbb R$ and each $i$. When $\gamma=1$, we get
	\begin{equation}
	Y_i^TY_i-\frac{(Y_i^T(\mathcal  AY_i)^{*3})^2}{||(\mathcal  AY_i)^{*3}||^2}\leq  Y_i^TY_i+2Y_i^T(\mathcal  AY_i)^{*3}+||(\mathcal  AY_i)^{*3}||^2\leq \sigma^2
	\end{equation}
	Therefore,
	\begin{equation}
	\lim\limits_{i\rightarrow\infty} \left(1-\frac{(Y_i^T(\mathcal  AY_i)^{*3})^2}{Y_i^TY_i||(\mathcal  AY_i)^{*3}||^2} \right)\leq \lim\limits_{n \rightarrow\infty} \frac{\sigma^2}{Y_i^TY_i}=0
	\end{equation}
	since $\{Y_i\}_i$ is unbounded. Let $X_i=\frac{Y_i}{||Y_i||}$, then $X_i\in \mathfrak U$, $(\mathcal  AX_i)^{*3}\neq \mathbf 0$ by the assumption $Y_i^T(\mathcal AY_i)^{*3}<0$ for all $i$, and then
	\begin{cor}\label{lem:dlimit}
		\begin{equation}
		\lim\limits_{i\rightarrow\infty} \frac{X_i^T(\mathcal  AX_i)^{*3}}{||(\mathcal  AX_i)^{*3}||}=\lim\limits_{i\rightarrow\infty} \frac{Y_i^T(\mathcal  AY_i)^{*3}}{||Y_i||\cdot||(\mathcal  AY_i)^{*3}||}=- 1 \label{eqn:keylim}
		\end{equation}	
	\end{cor}

	For each $i$, we decompose $X_i$ as $X_i=\alpha_i V_i+\sqrt{1-\alpha_i^2} W_i$ such that  $V_i\in \mathfrak U(\mathfrak A^\bot)$, $W_i\in \mathfrak U(Im(\mathcal A^T))$, and $\alpha_i\in [0,1]$. By assumption $Y_i^T(\mathcal AY_i)^{*3}<0$ for all $i$, we have
	\begin{cor}\label{cor:below1}
		$\alpha_i<1$ for all $i$.	
	\end{cor}
	
	Because $V_i,W_i\in \mathfrak U$ for all $i$, we can choose a subsequence $\{X_{i_k}\}_k$ of $\{X_i\}_i$ such that $\lim\limits_{k\rightarrow \infty} V_{i_k}=V_\infty\in \mathbb R^n$ and $\lim\limits_{k\rightarrow \infty} W_{i_k}=W_\infty\in \mathbb R^n$. As $ \mathfrak U(\mathfrak A^\bot)$ and $\mathfrak U(Im(\mathcal A^T))$ are closed, it must be $V_\infty \in \mathfrak U(\mathfrak A^{\bot})$ and $W_\infty \in \mathfrak U(Im(\mathcal A^T))$. Furthermore, we can choose a subsequence $\{X_{ik_j}\}_j$ of $\{X_{i_k}\}_k$ such that $\lim\limits_{j\rightarrow \infty} \alpha_{ik_j}=\alpha_\infty\in [0,1]$.  Without loss of generality, we assume $\{X_i\}_i$ has such properties of  $\{X_{ik_j}\}_j$, that is,
	\begin{enumerate}
		\item[(i)] $\lim\limits_{i\rightarrow\infty} \alpha_i=\alpha_\infty\in [0,1]$.
		
		\item[(ii)] $\lim\limits_{i\rightarrow\infty} V_i=V_\infty\in \mathfrak U(\mathfrak A^\bot)$ and $\lim\limits_{i\rightarrow\infty} W_i=W_\infty\in\mathfrak U(Im(\mathcal A^T))$
	\end{enumerate}
Then by the decomposition of $X_i$, each $Y_i$ can be uniquely decomposed as $Y_i=\alpha_i||Y_i||V_i+\sqrt{1-\alpha_i^2}||Y_i||W_i$. If $\{\sqrt{1-\alpha_i^2}||Y_i||W_i\}_i$   is bounded, then $\{\alpha_i||Y_i||V_i\}_i$ and so $\{\mathcal C(Y_i)=\alpha_i||Y_i|| V_i+\sqrt{1-\alpha_i^2}||Y_i||W_i+(\sqrt{1-\alpha_i^2}||Y_i||)^3(\mathcal A W_i)^{*3}\}_i$ must be  unbounded since $\{Y_i\}_i$ is unbounded. Therefore,  $\{\sqrt{1-\alpha_i^2}||Y_i||W_i\}_i$ and so $\{\sqrt{1-\alpha_i^2}||Y_i||\}_i$  must be unbounded. As a result, 
\begin{cor}
$\alpha_i>0$ for infinitely many $i$.
\end{cor} Without loss of generality, we assume $\alpha_i>0$ for all $i$.  Note that $\lim\limits_{i\rightarrow \infty} \frac{\sqrt{1-\alpha_i^2}||Y_i||\}_iW_i}{||\sqrt{1-\alpha_i^2}||Y_i||W_i||}=W_\infty \in \mathfrak U(Im(\mathcal A^T))\subset Im(\mathcal A^T)$ and is nonzero. Thus, $\mathcal A W_\infty\neq \mathbf 0$. Then we have
\begin{eqnarray}
\mathbf 0 & =& \lim\limits_{i\rightarrow\infty} \mathcal A\left(\frac{\mathcal C(Y_i)}{\sqrt{1-\alpha_i^2}||Y_i||}\right) \\
& =& \mathcal A W_\infty + \lim\limits_{i\rightarrow\infty}\left(\sqrt{1-\alpha_i^2}||Y_i||\right)^2\mathcal A(\mathcal A W_i)^{*3}
\end{eqnarray}
In consequence, $W_\infty\in Im(\mathcal A^T)$ and $\mathcal A(\mathcal A W_\infty)^{*3}=\lim\limits_{i\rightarrow\infty}\mathcal A(\mathcal A W_i)^{*3}=\mathbf 0$ since $\{\sqrt{1-\alpha_i^2}||Y_i||\}_i$ is unbounded.  Let's take $W:=W_\infty$, then $W\in Im(\mathcal A^T)$, $\mathcal A(\mathcal A W)^{*3}=\mathbf 0$, and $\mathcal A W\neq \mathbf 0$.

Let $U_i:=\sqrt{1-\alpha_i^2}||Y_i|| \mathcal A W_i\in  Im(\mathcal A\mathcal A^T)$, since $W_i\in Im(\mathcal A^T)$.	Let $\delta_1:=\inf \{||\mathcal A X||\mid X\in Im (\mathcal A^T) \& ||X||=1\}$, then $\delta_1>0$ by Propositions \ref{prop:property} and \ref{prop:unitclosed}, and $\delta_1 ||X||\leq ||\mathcal A X||$ for all $X\in Im(\mathcal A^T)$. Thus, $\{ U_i\in Im(\mathcal A\mathcal A^T)\}_i$ is unbounded since $\{\sqrt{1-\alpha_i^2}||Y_i||  W_i\}$ is unbounded. Notice that
\begin{eqnarray}
 U_i +\mathcal A U_i^{*3} & =& \mathcal A (Y_i + (\mathcal A Y_i)^{*3}) \\
 & = & \mathcal A \cdot\mathcal C (Y_i)
\end{eqnarray}
By the assumption on $\{Y_i\}_i$, $\{\mathcal C (Y_i)\}_i$ is bounded. Therefore, $\{ U_i +\mathcal A U_i^{*3}\}_i$ is bounded. In consequence, $\hat{\mathcal C}(X)$ is not proper on $Im(\mathcal A\mathcal A^T)$.

This complete the proof of Lemma \ref{lem:nonprop}.
\end{IEEEproof}	

Note that this result is also independently obtained by Tuyen in \cite{Tuyen2020}. In fact, the inverse of (I) is also true, please refer to \cite{Tuyen2020}. Now, we are ready to prove Theorem \ref{thm:D-proper} by contradiction.

Assume there is some unbounded sequence $\{U_i \}_i\subset Im(\mathcal A \mathcal A^T)$ such that $\{U_i+\mathcal A U_i^{*3}\}_i$ is bounded. Similar to before, there is some positive real number $\beta>0$ such that
\begin{eqnarray}
\beta & \geq & ||U_i+\mathcal A U_i^{*3}||^2 \\
 &=& ||U_i||^2 +2U_i^T\mathcal A U_i^{*3} +||\mathcal A U_i^{*3}||^2 \\
 &\geq & ||U_i||^2 -\frac{(U_i^T\mathcal A U_i^{*3})^2}{||\mathcal A U_i^{*3}||^2} 
\end{eqnarray}
Herein, $||\mathcal A U_i^{*3}||>0$ for almost all $i$ since $\beta<\infty$.
In consequence, we have
\begin{equation}
\lim\limits_{i\rightarrow\infty} \frac{U_i^T\mathcal A U_i^{*3}}{||U_i||\cdot||\mathcal A U_i^{*3}||}=-1
\end{equation}
As $\{\frac{U_i}{||U_i||}\}_i\subset \mathfrak U$, we can assume it has a limit. Let $U_\infty:=\lim\limits_{i\rightarrow\infty}\frac{U_i}{||U_i||}\in \mathfrak U$. It is evident that $Im(\mathcal A \mathcal A^T)$ is a closed set. So $U_\infty\in Im(\mathcal A \mathcal A^T)$. Similarly, we consider a limit of $\{\frac{\mathcal A U_i^{*3}}{||\mathcal A U_i^{*3}||}\}_i$. Without loss of generality, we can assume $\frac{\mathcal A U_\infty^{*3}}{||\mathcal A U_\infty^{*3}||}=\lim\limits_{i\rightarrow\infty}\frac{\mathcal A U_i^{*3}}{||\mathcal A U_i^{*3}||}$. Then we have
\begin{equation}\label{eqn:limf1}
\frac{U_\infty^T \mathcal AU_\infty^{*3} }{||\mathcal AU_\infty^{*3}||}=\lim\limits_{i\rightarrow\infty} \frac{U_i^T\mathcal A U_i^{*3}}{||U_i||\cdot||\mathcal A U_i^{*3}||}=-1
\end{equation}
Thus, there is some number $\gamma>0$ such that
\begin{equation}\label{eqn:zero}
 U_\infty+\gamma \mathcal AU_\infty^{*3}=\mathbf 0
\end{equation}
By $U_\infty\in Im(\mathcal A \mathcal A^T)$, we take a $Z_\infty\in Im(\mathcal A^T)$ such that $U_\infty=\mathcal A Z_\infty$. By equation (\ref{eqn:zero})
\begin{eqnarray}
\mathbf 0 &=& U_\infty+\gamma \mathcal AU_\infty^{*3} \\
 &= & (\mathcal I + \gamma \mathcal A (U_\infty)^{\Delta 2}) U_\infty \\
 &= & (\mathcal I + \gamma \mathcal A (\mathcal A Z_\infty)^{\Delta 2}) \mathcal A Z_\infty
\end{eqnarray}
This is impossible, since $(\mathcal I + \gamma \mathcal A (\mathcal A Z_\infty)^{\Delta 2}) $ is invertible  by Proposition \ref{prop:d-mapclas} and $U_\infty=\mathcal A Z_\infty \in \mathfrak U$ is not zero. Therefore, for a D-map $\mathcal D(X):=X+(\mathcal A X)^{*3}$, the map $\hat{\mathcal D}(X):=X+\mathcal A X^{*3}$ must be proper on $Im(\mathcal A \mathcal A^T)$. In consequence, each D-map must be proper.

\section{Discussion and Related Work}\label{sec:dis-rel}
Recall the history of $\mathcal{JC}$, there are many excellent works. First of all,  Keller proposed  $\mathcal{JC}$ for $\mathbb K^n=\mathbb C^2$ in 1939 \cite{Keller1939}. And Abhyankar gave its modern style in his lectures \cite{Abhyankar1977,Abhyankar1974Inversion}.
Let $(\mathcal{JC})^k$ denote $\mathcal{JC}(\mathbb K , n)$ in which the polynomial degrees are not greater than $k$.  A remarkable progress is Wang's result \cite{Wang1980degree2} that  $(\mathcal{JC})^2$ is true for arbitrary field of characteristic $\neq 2$. This led to the studies on the reduction of $\mathcal{JC}$ to specific $(\mathcal{JC})^k$ for some small integer $k$. In the line of this, Yag{\v{z}}ev \cite{Yagzhev1981} and independently Bass et al. \cite{Bass1982formalinverse} showed that  $(\mathcal{JC})^3$ implies $\mathcal{JC}$. Furthermore, Dru{\.{z}}kowski \cite{druzkowski1983} showed that it suffices to show D-maps being automorphism for  $\mathcal{JC}$. Anyway, it cannot make a reduction of $\mathcal{JC}$ to  $(\mathcal{JC})^2$ \cite{Campbell1973} on  the field $\mathbb C$. That is, we must solve $(\mathcal{JC})^3$ to the end of $\mathcal{JC}$. Anyway, many studies \cite{Ludwik2005,Essen2005HC} showed that it suffices to prove specific structure D-maps for  $\mathcal{JC}$. Furthermore, $\mathcal{JC}$  also appeared to be connected to questions in noncommutative algebra, for example, $\mathcal{JC}$ is equivalent to the  Dixmier Conjecture which asserts that each endomorphism of the Weyl algebra is surjective (hence an automorphism) \cite{tsuchimoto2005DC,belovkanel2007DC,Adjamagbo2007DJC}. $\mathcal{JC}$ is proved equivalent to various conjectures, such as,   Kernel Conjecture  \cite{denessen2001review},  Hessian Conjecture \cite{meng2006HC,Essen2005HC}, Eulerian Conjecture \cite{Adjamagbo1991EC}, etc. There are many  partial results of $\mathcal{JC}$ on special categories of polynomials, for instance, the ``non-negative coefficients" D-map \cite{Druzkowski1997PositiveCoefficients}, D-maps in low dimension space \cite{cheng20024D,deondt20053D}, a special class of D-maps in dimension $9$ \cite{Yan2011}, tame automorphisms \cite{smith1989stablytame,drensky2000tame,shestakov2004tame}, etc, refer to surveys \cite{Bass1982formalinverse,Druzkowski1997ReviewJC,Essenbook2000,denessen2001review,Wright2005Review} for more related  results. There were some studies about $\mathcal{JC}$ for fixed number of variables, even for  $2$-variables, such as,  $\mathcal{JC}$ \cite{Moh1983deg100} for $2$-variables $\leq 100$-degree, sufficient conditions via polynomial flows in \cite{bass1985polynomialflows} for $\mathcal{JC}(\mathbb R, 2)$, a Hamiltonian flows approach in \cite{campbell1997jacobian} for $\mathcal{JC}(\mathbb C, 2)$.

$\mathcal{JC}$ must depend on Jacobian Condition, polynomial type, and the number field of characteristic $0$. For the Keller maps over number fields of characteristic $0$, injectivity always implies subjectivity \cite{Bailynicki1962Injective,Winiarski1979Inverse,Yagzhev1981,Bass1982formalinverse,Kurdyka1988Surjectivity}.
But, this property completely fails for the nonpolynomial maps, already for $n=2$. There is a counterexample in \cite{Bass1982formalinverse}: $ F_1(X)=e^{x_1}, F_{2}=x_2e^{-x_1}$ whose Jacobian is $1$, but $\mathcal F (\mathbb C^2)$ excludes exactly the axis $x_1=0$. That is, injectivity does not mean surjectivity for a generic analytic map. Even for rational maps, there are counterexamples in \cite{druzkowski1995SurveyCounterexamples}. As to the zero characteristic condition, there are counterexamples in  \cite{Abhyankar1974Inversion,Bass1982formalinverse} for $\mathcal{JC}$ of characteristic $>0$. The Jacobian Condition also cannot be relaxed. A generalization of $\mathcal {JC}$ is the real Jacobian problem \cite{Randall1983RJCP} (also called strong Jacobian Conjecture in \cite{Pinchuk1994}),  that is, whether a polynomial mapping $\mathcal F:\mathbb R^2\mapsto \mathbb R^2$ with a nonvanishing Jacobian determinant is an automorphism. The strong Jacobian Conjecture  has a negative answer \cite{Pinchuk1994}. Pinchuk presented
a beautiful example of a non-injective polynomial mapping $\mathcal F (x_1,x_2)$ of $\mathbb R^2$ into itself, of degree $(x_1,x_2)=(10, 40)$,  whose Jacobian determinant is everywhere positive on $\mathbb R^2$.  Therefore, Polynomial Automorphism Theorem (PAT for short) is the best in all of what we can get.

From PAT, it immediately gets that the Roll Theorem is true for polynomial functions over any algebraically closed field. Let $\mathbb K $ be an algebraically closed field of characteristic zero, $\mathcal P\in \mathbb K [X]^n$. Then  the determinant $det(\mathbf J_{\mathcal P}(X))$  is a polynomial and  must have an $X_0\in\mathbb K^n$ such that $det(\mathbf J_{\mathcal P}(X_0))=0$ if $det(\mathbf J_{\mathcal P}(X))$ is not a constant.
By Theorem \ref{thm:Conjecture}, using contradiction argument we can obtain
\begin{HDRthm}\label{thm:HDRoll}
	If there are $X\neq Y\in \mathbb K ^n$ such that $\mathcal P(X)= \mathcal P(Y)$ then there is some $X_0$ such that $det(\mathbf J_{\mathcal P}(X_0))=0$.
\end{HDRthm}
In this theorem, the requirement ``algebraically closed" is necessary. Otherwise, the Pinchuk's counterexample will be a counterexample over $\mathbb R^2$.

Another immediate result of PAT is  that the inverse flows  are actually polynomials in $X$ and $t$. So, high order Lie derivatives vanish at some stages. That is, the Lie derivatives are locally nilpotent or finite \cite{VANDENESSEN1992Locally,Arno1994Locally}. Therefore, it actually gives a termination criterion for computing inverse polynomial through Lie derivatives  \cite{nousiainen1983automorphisms}.

When a polynomial map is an automorphism, there are several different approaches to the inversion formulas. An early one is the  Abhyankar-Gurjar inversion formula \cite{Abhyankar1974Inversion}. In   \cite{Bass1982formalinverse} Bass et al.  presented a formal expansion for the inverse.  Nousiainen and Sweedler \cite{nousiainen1983automorphisms} provided an inversion formula through Lie derivatives. For specific polynomials, Wright \cite{wright1987formal} and respectively Zhao \cite{Zhaorecurrentinversion,Zhao2005Inversion} gave advanced inversion formals. Anyway, the degree of inverse polynomial is bounded by  $deg(F^{-1})\leq deg(F)^{n-1}$ \cite{Bass1982formalinverse,Rusek1984Polynomial,furter1998inversdegree}.

Besides Jacobian Condition, van den Essen \cite{van1990aGrobnerInversion}  using Gr\"obner base  gave an algebraic criterion for the invertibility of polynomial maps.
\begin{Ethm}[\cite{van1990aGrobnerInversion}]\label{thm:Inverse-Grobner}
	Any  map $\mathbf F=(F_1(X),\dots,F_{n}(X))\in \mathbb K [X]$ on arbitrary field $\mathbb K $ is an automorphism iff there are polynomials $G_1(Y), \dots$, $G_n(Y)\in \mathbb K [Y]$ such that $Y_1-F_1(X), \dots, Y_n-F_n(X)$ and $X_1-G_1(Y), \dots, X_n-G_n(Y)$ generate the same polynomial idea in $\mathbb K [X, Y]$.
\end{Ethm}

This criterion is not limited to the characteristic zero cases but holds in all characteristics. At the same time Theorem \ref{thm:Inverse-Grobner} also provides an algorithm to decide if a polynomial map has an inverse and compute the inverse if it exists. The theory of Gr\"obner bases for polynomial ideals \cite{Cox1997Ideals} is the foundation of the Essen Theorem. In contrast with this, PAT is an analytic criterion for the global invertibility of polynomial maps. In particular, PAT can be efficiently implemented for sparse polynomial maps, by testing the Jacobian Condition through random inputs.

\section{Concluding Remarks}\label{sec:con}
In this study, we gave an affirmative answer to Jacobian Conjecture.  Based on D-map, our proof used algebraic property like nilpotency property. To study the direction changing tendency of the unbounded sequence $\{Y_i\}_i$, we used an optimization method to obtain the key limit equations (\ref{eqn:keylim}), (\ref{eqn:limf1}) and (\ref{eqn:zero}). Jacobian Conjecture is an algebraic geometry problem. It is no surprise to use algebraic methods in the proof, but the optimization method is really an extra auxiliary. So, this proof is said as an optimization-based method. Our proof takes much advantage of D-maps' nice algebraic  properties. Although Yag{\v{z}}ev map is a little extension of D-map, it has no such good properties. At least, so far our proof does not work on Yag{\v{z}}ev map for which we cannot get a necessary condition to nonproperness like Lemma \ref{lem:nonprop}.  From the proof, we can clearly see how and why Jacobian Condition makes  D-map being proper and so injective, like the proof for  quadratic polynomial maps. However, we should note that an injective polynomial map in $\mathbb R[X]^n$ is not necessarily an automorphism. The automorphism property of  D-maps in $\mathbb R[X]^n$ is not a direct result of their injectivity but  derived from  a series of deductions involved with Lefschetz Principle, Dru{\.{z}}kowski's reduction, automorphism property of  D-maps in $\mathbb C[X]^n$, etc.

Given a polynomial map $\mathcal P(X):=(\mathcal P_1(X),\dots, \mathcal P_n(X))\in \mathbb K[X]^n$, $\mathcal P(X)$ is an automorphism from $\mathbb K^n$ onto $\mathbb K^n$ iff the induced endomorphism $\mathcal R_{\mathcal P}:  \mathbb K[X]^n \mapsto \mathbb K[X]^n$ by $\mathcal R_{\mathcal P}(X_i)=\mathcal P_i(X)$ for $i=1,\dots, n$ is an automorphism of the ring $\mathbb K[X]^n$. If we consider the derivatives as algebraic operators, then the Jacobian Conjecture is purely an algebraic problem. So a purely algebraic proof is really an interesting thing. It is already known that an analytic map satisfying Jacobian Condition is not necessarily a diffeomorphism. In fact, even for the rational maps, the Jacobian Condition is not a sufficient condition for this type of map being diffeomorphisms. For analytic maps over Euclidean space, Hadamard's Diffeomorphism Theorem has provided a nice criterion for diffeomorphism. However, this is not a computable approach like Jacobian Condition. To the best of our knowledge, so far there is no computable method to directly check the properness of a map. In practice, we may need to computably determine whether a given concrete map is a diffeomorphism. In such context, Hadamard's Diffeomorphism Theorem helps a little and the Jacobian Condition is not a correct criterion for nonpolymomial maps. In the physical world, we are usually concerned with elementary maps which are composed of elementary expressions like $e^X, X^p, \sin X, \cos X$, etc. So for analytic maps, it is natural to ask whether there is some computable diffeomorphism criterion. Another basic question about polynomial automorphisms is how many of them, or what is the ratio of polynomial automorphisms to all polynomials, or what is their distribution. By Weierstrass Approximation Theorem, each continuous real function on some closed interval can be uniformly approximated by polynomials. In comparison with this, it is an interesting problem whether the automorphism polynomials are dense in the set of all analytic diffeomorphisms. The study of the injectivity has given rise to several surprising results and interesting relations in various directions and from different perspectives. In this study, we proved Jacobian Conjecture by showing the injectivity of D-maps. In fact, the injectivity itself has received attention from not only the mathematical field but also by the economic field \cite{Samuelson1953Prices}, physical field \cite{abdesselam2003QFT}, and chemical field \cite{muller2019bijectivity}. This study only verifies the injectivity of a special class of polynomial maps in $\mathbb R[X]^n$. The proof heavily depends on the Jacobian Condition and the form of D-map. It has been proved that the Samuelson Conjecture in \cite{Samuelson1953Prices} is true for any polynomial map \cite{denessen1992Samuelsonpolynomial}  and arbitrary rational map \cite{campbell1994rational}, but fails for generic analytic maps \cite{Gale1965GlobalUnivalence}. For a general differentiable map $\mathbf F$ on  $\mathbb R^n$, Chamberland Conjecture \cite{chamberland1998conjecture}  asserts if all the eigenvalues of $\mathbf J_{\mathbf F}(X)$ are bounded away from zero then $\mathbf F$ must be injective. This conjecture is still open.

\section{Acknowledgement}
The author would like to thank Tuyen Trung Truong for his counterexample \cite{Tuyenexm} to a previous aproach. The author also would like to thank Chengling Fang for her help to find out the exact error in the previously invalid proof. This work was partially supported by NSF of China (No. 61672488), CAS Youth Innovation Promotion Association (No. 2015315), 
National Key R$\&$D Program of China (No. 2018YFC0116704).

\bibliographystyle{IEEEtran}

\begin{thebibliography}{10}
	\providecommand{\url}[1]{#1}
	\csname url@samestyle\endcsname
	\providecommand{\newblock}{\relax}
	\providecommand{\bibinfo}[2]{#2}
	\providecommand{\BIBentrySTDinterwordspacing}{\spaceskip=0pt\relax}
	\providecommand{\BIBentryALTinterwordstretchfactor}{4}
	\providecommand{\BIBentryALTinterwordspacing}{\spaceskip=\fontdimen2\font plus
		\BIBentryALTinterwordstretchfactor\fontdimen3\font minus
		\fontdimen4\font\relax}
	\providecommand{\BIBforeignlanguage}[2]{{%
			\expandafter\ifx\csname l@#1\endcsname\relax
			\typeout{** WARNING: IEEEtran.bst: No hyphenation pattern has been}%
			\typeout{** loaded for the language `#1'. Using the pattern for}%
			\typeout{** the default language instead.}%
			\else
			\language=\csname l@#1\endcsname
			\fi
			#2}}
	\providecommand{\BIBdecl}{\relax}
	\BIBdecl
	

\bibitem{abdesselam2003QFT}
\bgroup\scshape{}A.~Abdesselam\egroup{}, The {J}acobian conjecture as a problem
of perturbative quantum field theory,  \emph{Annales Henri Poincaré}
\textbf{4} (2003), 199--215.

\bibitem{Abhyankar1977}
\bgroup\scshape{}S.~S. Abhyankar\egroup{}, Lectures on expansion techniques in
algebraic geometry. with notes by balwant singh,  \emph{Tata Institute of
	Fundamental Research Bombay} (1977).

\bibitem{Abhyankar1974Inversion}
\bgroup\scshape{}S.~S. Abhyankar\egroup{}, Lectures in algebraic geometry,
\emph{Notes by Chris Christensen,Purdue University} (1974).

\bibitem{Adjamagbo1991EC}
\bgroup\scshape{}K.~Adjamagbo\egroup{} and \bgroup\scshape{}A.~van~den
Essen\egroup{}, Eulerian systems of partial differential equations and the
{J}acobian conjecture,  \emph{Journal of Pure and Applied Algebra}
\textbf{74} (1991), 1--15.

\bibitem{Adjamagbo2007DJC}
\bgroup\scshape{}P.~K. Adjamagbo\egroup{} and \bgroup\scshape{}A.~V.~D.
Essen\egroup{}, A proof of the equivalence of the {D}ixmier, {J}acobian and
{P}oisson conjectures,  \emph{Acta Mathematica Vietnamica} \textbf{32}
(2007), 205--14.

\bibitem{Bailynicki1962Injective}
\bgroup\scshape{}A.~Bailynicki-Birula\egroup{} and
\bgroup\scshape{}M.~Rosenlicht\egroup{}, Injective morphisms of real
algebraic varieties,  \emph{Proceedings of the American Mathematical Society}
\textbf{13} (1962), 200--203.

\bibitem{Bass1982formalinverse}
\bgroup\scshape{}H.~Bass\egroup{}, \bgroup\scshape{}E.~H. Connell\egroup{}, and
\bgroup\scshape{}D.~Wright\egroup{}, The {J}acobian conjecture: Reduction of
degree and formal expansion of the inverse,  \emph{Bulletin of the American
	Mathematical Society} \textbf{7} (1982), 287--330.

\bibitem{bass1985polynomialflows}
\bgroup\scshape{}H.~Bass\egroup{} and \bgroup\scshape{}G.~H. Meisters\egroup{},
Polynomial flows in the plane,  \emph{Advances in Mathematics} \textbf{55}
(1985), 173--208.

\bibitem{belovkanel2007DC}
\bgroup\scshape{}A.~Belovkanel\egroup{} and
\bgroup\scshape{}M.~Kontsevich\egroup{}, The {J}acobian conjecture is stably
equivalent to the {D}ixmier conjecture,  \emph{Moscow Mathematical Journal}
\textbf{7} (2007), 209--218.

\bibitem{deondt20053D}
\bgroup\scshape{}M.~D. Bondt\egroup{} and \bgroup\scshape{}A.~van~den
Essen\egroup{}, The {J}acobian conjecture: Linear triangularization for
homogeneous polynomial maps in dimension three,  \emph{Journal of Algebra}
\textbf{294} (2005), 294--306.

\bibitem{Boyd2004Optimization}
\bgroup\scshape{}S.~Boyd\egroup{} and
\bgroup\scshape{}L.~Vandenberghe\egroup{}, \emph{Convex Optimization},
Cambridge University Press, 2004.

\bibitem{campbell1994rational}
\bgroup\scshape{}L.~A. Campbell\egroup{}, Rational {S}amuelson maps are
univalent,  \emph{Journal of Pure and Applied Algebra} \textbf{92} (1994),
227--240.

\bibitem{Campbell1973}
\bgroup\scshape{}L.~A. Campbell\egroup{}, A condition for a polynomial map to
be invertible,  \emph{Mathematische Annalen} \textbf{205} (1973), 243--248.

\bibitem{campbell1997jacobian}
\bgroup\scshape{}L.~A. Campbell\egroup{}, Jacobian pairs and {H}amiltonian
flows,  \emph{Journal of Pure and Applied Algebra} \textbf{115} (1997),
15--26.

\bibitem{chamberland1998conjecture}
\bgroup\scshape{}M.~Chamberland\egroup{} and
\bgroup\scshape{}G.~Meisters\egroup{}, A mountain pass to the {J}acobian
conjecture.,  \emph{Canadian Mathematical Bulletin} \textbf{41} (1998),
442--451.

\bibitem{cheng20024D}
\bgroup\scshape{}C.~C. Cheng\egroup{}, Power linear {K}eller maps of dimension
four,  \emph{Journal of Pure and Applied Algebra} \textbf{169} (2002),
153--158.

\bibitem{Cox1997Ideals}
\bgroup\scshape{}D.~A. Cox\egroup{}, \bgroup\scshape{}J.~Little\egroup{}, and
\bgroup\scshape{}D.~O'Shea\egroup{}, \emph{Ideals, Varieties, and
	Algorithms}, Springer, 1997.

\bibitem{denessen1992Samuelsonpolynomial}
\bgroup\scshape{}A.~V. Den~Essen\egroup{} and
\bgroup\scshape{}T.~Parthasarathy\egroup{}, Polynomial maps and a conjecture
of {S}amuelson,  \emph{Linear Algebra and its Applications} \textbf{177}
(1992), 191--195.

\bibitem{drensky2000tame}
\bgroup\scshape{}V.~Drensky\egroup{}, \bgroup\scshape{}A.~van~den
Essen\egroup{}, and \bgroup\scshape{}D.~Stefanov\egroup{}, New stably tame
automorphisms of polynomial algebras,  \emph{Journal of Algebra} \textbf{226}
(2000), 629--638.

\bibitem{Druzkowski1997ReviewJC}
\bgroup\scshape{}L.~M. Dru{\.{z}}kowski\egroup{}, Partial results and
equivalent formulations of the {J}acobian conjecture,  \emph{Rendiconti Del
	Seminario Matematico} (1997), 275--282.

\bibitem{druzkowski1983}
\bgroup\scshape{}L.~M. Dru{\.{z}}kowski\egroup{}, An effective approach to
{K}eller's {J}acobian conjecture,  \emph{Mathematische Annalen} \textbf{264}
(1983), 303--313.

\bibitem{druzkowski1995SurveyCounterexamples}
\bgroup\scshape{}L.~M. Dru{\.{z}}kowski\egroup{}, The {J}acobian conjecture: survey of some results,  \emph{Banach
	Center Publications} \textbf{31} (1995), 163--171.

\bibitem{Druzkowski1997PositiveCoefficients}
\bgroup\scshape{}L.~M. Dru{\.{z}}kowski\egroup{}, The {J}acobian conjecture in case of ``non-negative coefficients",
\emph{Annales Polonici Mathematici} \textbf{66} (1997), 67--75.

\bibitem{Ludwik2005}
\bgroup\scshape{}L.~M. Dru{\.{z}}kowski\egroup{}, The {J}acobian conjecture:
symmetric reduction and solution in the symmetric cubic linear case,
\emph{Annales Polonici Mathematici} \textbf{87} (2005), 83--92.

\bibitem{van1990aGrobnerInversion}
\bgroup\scshape{}A.~van~den Essen\egroup{}, A criterion to decide if a
polynomial map is invertible and to compute the inverse,
\emph{Communications in Algebra} \textbf{18} (1990), 3183--3186.

\bibitem{VANDENESSEN1992Locally}
\bgroup\scshape{}A.~van~den Essen\egroup{}, Locally finite and locally nilpotent derivations with applications to
polynomial flows and polynomial morphisms,  \emph{Proceedings of the American
	Mathematical Society} \textbf{116} (1992), 861--871.

\bibitem{Arno1994Locally}
\bgroup\scshape{}A.~van~den Essen\egroup{}, Locally finite and locally nilpotent derivations with applications to
polynomial flows, morphisms and $\mathcal{Y}_a$-actions. {II}.,
\emph{Proceedings of the American Mathematical Society} \textbf{121} (1994),
667--678.

\bibitem{Essenbook2000}
\bgroup\scshape{}A.~van~den Essen\egroup{}, Polynomial automorphisms and the {J}acobian conjecture,  in
\emph{Progress in Mathematics Vol 190}, Boston; Berlin: Birkh{\"a}user, 2000.

\bibitem{denessen2001review}
\bgroup\scshape{}A.~van~den Essen\egroup{}, Polynomial automorphisms, and the {J}acobian conjecture,  \emph{The
	Mathematical Gazette} \textbf{85} (2001), 572.

\bibitem{Essen2005HC}
\bgroup\scshape{}A.~van~den Essen\egroup{}, A reduction of the {J}acobian conjecture to the symmetric case,
\emph{Proceedings of the American Mathematical Society} \textbf{133} (2005),
2201--2205.

\bibitem{furter1998inversdegree}
\bgroup\scshape{}J.~Furter\egroup{}, On the degree of the inverse of an
automorphism of the affine space,  \emph{Journal of Pure and Applied Algebra}
\textbf{130} (1998), 277--292.

\bibitem{Gale1965GlobalUnivalence}
\bgroup\scshape{}D.~Gale\egroup{} and \bgroup\scshape{}H.~Nikaido\egroup{}, The
{J}acobian matrix and global univalence of mapping,  \emph{Mathematische
	Annalen} \textbf{159} (1965), 81--93.

\bibitem{Gordon1972On}
\bgroup\scshape{}W.~B. Gordon\egroup{}, On the diffeomorphisms of {E}uclidean
space,  \emph{American Mathematical Monthly} \textbf{79} (1972), 755--759.

\bibitem{Hadamard1906}
\bgroup\scshape{}J.~Hadamard\egroup{}, Sur les transformations ponctuelles,
\emph{Bull.soc.math.france} (1906), 71--84.

\bibitem{Keller1939}
\bgroup\scshape{}O.~Keller\egroup{}, Ganze {C}remona-transformationen,
\emph{Monatshefte f\"{u}r Mathematik} \textbf{47} (1939), 299--306.

\bibitem{Kurdyka1988Surjectivity}
\bgroup\scshape{}K.~Kurdyka\egroup{} and \bgroup\scshape{}K.~Rusek\egroup{},
Surjectivity of certain injective semialgebraic transformations of
$\mathbb{R}^n$,  \emph{Mathematische Zeitschrift} \textbf{200} (1988),
141--148.

\bibitem{meng2006HC}
\bgroup\scshape{}G.~Meng\egroup{}, Legendre transform, {H}essian conjecture and
tree formula,  \emph{Applied Mathematics Letters} \textbf{19} (2006),
503--510.

\bibitem{Moh1983deg100}
\bgroup\scshape{}T.~T. Moh\egroup{}, On the {J}acobian conjecture and the
configurations of roots,  \emph{Journal F{\"u}r Die Reine Und Angewandte
	Mathematik} \textbf{1983} (1983), 140--213.

\bibitem{muller2019bijectivity}
\bgroup\scshape{}S.~Muller\egroup{}, \bgroup\scshape{}J.~Hofbauer\egroup{}, and
\bgroup\scshape{}G.~Regensburger\egroup{}, On the bijectivity of families of
exponential/generalized polynomial maps,  \emph{SIAM Journal on Applied
	Algebra and Geometry} \textbf{3} (2019), 412--438.

\bibitem{nousiainen1983automorphisms}
\bgroup\scshape{}P.~Nousiainen\egroup{} and
\bgroup\scshape{}M.~Sweedler\egroup{}, Automorphisms of polynomial and power
series rings,  \emph{Journal of Pure and Applied Algebra} \textbf{29} (1983),
93--97.

\bibitem{Pinchuk1994}
\bgroup\scshape{}S.~Pinchuk\egroup{}, A counterexamle to the strong real
{J}acobian conjecture,  \emph{Mathematische Zeitschrift} \textbf{217} (1994),
1--4.

\bibitem{Randall1983RJCP}
\bgroup\scshape{}J.~D. Randall\egroup{}, The real {J}acobian problem, 1983,
pp.~411--414.

\bibitem{RudinInjective1995}
\bgroup\scshape{}W.~Rudin\egroup{}, Injective polynomial maps are
automorphisms,  \emph{American Mathematical Monthly} \textbf{102} (1995),
540--543.

\bibitem{Rusek1984Polynomial}
\bgroup\scshape{}K.~Rusek\egroup{} and \bgroup\scshape{}T.~Winiarski\egroup{},
Polynomial automorphisms of $\mathbb{C}^n$,  \textbf{661} (1984).

\bibitem{Samuelson1953Prices}
\bgroup\scshape{}P.~A. Samuelson\egroup{}, Prices of factors and good in
general equilibrium,  \emph{Review of Economic Studies} \textbf{21} (1953),
1--20.

\bibitem{shestakov2004tame}
\bgroup\scshape{}I.~Shestakov\egroup{} and
\bgroup\scshape{}U.~Umirbaev\egroup{}, The tame and the wild automorphisms of
polynomial rings in three variables,  \emph{Journal of the American
	Mathematical Society} \textbf{17} (2004), 197--227.

\bibitem{smith1989stablytame}
\bgroup\scshape{}M.~K. Smith\egroup{}, Stably tame automorphisms,
\emph{Journal of Pure and Applied Algebra} \textbf{58} (1989), 209--212.

\bibitem{Tuyen2020}
\bgroup\scshape{}T.~T. Truong\egroup{}, Some observations on the properness of identity plus linear powers, arXiv:2004.03309, 2020.

\bibitem{Tuyenexm}
\bgroup\scshape{}T.~T. Truong\egroup{}, Some observations on the properness of identity plus linear powers:
part 2,  arXiv:2005.02260, 2020.

\bibitem{tsuchimoto2005DC}
\bgroup\scshape{}Y.~Tsuchimoto\egroup{}, Endomorphisms of {W}eyl algebra and
$p$-curvatures,  \emph{Osaka Journal of Mathematics} \textbf{42} (2005),
435--452.

\bibitem{Wang1980degree2}
\bgroup\scshape{}S.~S. Wang\egroup{}, A {J}acobian criterion for separability,
\emph{Journal of Algebra} \textbf{65} (1980), 453--494.

\bibitem{Winiarski1979Inverse}
\bgroup\scshape{}T.~Winiarski\egroup{}, Inverse of polynomial automorphisms of
$\mathbb{C}^n$,  \emph{Journal of the Institute of Image Information $\&$
	Television Engineers} \textbf{27} (1979), 9--11.

\bibitem{wright1987formal}
\bgroup\scshape{}D.~Wright\egroup{}, Formal inverse expansion and the
{J}acobian conjecture,  \emph{Journal of Pure and Applied Algebra}
\textbf{48} (1987), 199--219.

\bibitem{Wright2005Review}
\bgroup\scshape{}D.~Wright\egroup{}, The {J}acobian conjecture: ideal membership questions and recent
advances,  (2005), 261--276.

\bibitem{Yagzhev1981}
\bgroup\scshape{}A.~V. Yag{\v{z}}ev\egroup{}, Keller's problem,  \emph{Siberian
	Mathematical Journal} \textbf{21} (1980), 747--754.

\bibitem{Yan2011}
\bgroup\scshape{}D.~Yan\egroup{}, A note on the {J}acobian conjecture,
\emph{Linear Algebra and its Applications} \textbf{435} (2011), 2110--2113.

\bibitem{Zhaorecurrentinversion}
\bgroup\scshape{}W.~Zhao\egroup{}, Recurrent inversion formulas,  in \emph{Some
	Properties and Open Problems of {H}essian Nilpotent polynomials, In
	preparation. Department of Mathematics, Illinois State University, Normal, IL
	61790-4520. E-mail: wzhao@ilstu.edu}.

\bibitem{Zhao2005Inversion}
\bgroup\scshape{}W.~Zhao\egroup{}, Inversion problem, {L}egendre transform and inviscid {B}urgers'
equations,  \emph{Journal of Pure and Applied Algebra} \textbf{199} (2005),
299--317.

	
\end{thebibliography}

\end{document}